\documentclass[12pt]{article}
 \usepackage{amssymb}
\newtheorem{theorem}{Theorem}

\newtheorem{proposition}[theorem]{Proposition}
\newtheorem{lemma}[theorem]{Lemma}

\newtheorem{remark}[theorem]{Remark}

\newtheorem{example}[theorem]{Example}

\newtheorem{definition}[theorem]{Definition}

\newcommand{\R}{\mathbb{R}}
\newcommand{\Ee}{\mathbb {E}}
\newcommand{\Q}{\mathbb{Q}}

\newcommand{\Sf}{\mathbb{S}}

\newcommand{\Le}{\mathbb{L}}

\newcommand{\Hy}{\mathbb{H}}

\newcommand{\spa}{\mbox{span}}

\newcommand{\po}{{\hspace*{-1ex}}{\bf .  }}

\def\<{\langle}

\def\>{\rangle}
\def\a{\alpha}

\def\e{\epsilon}
\def\d{\partial}
\def\bea{\begin{eqnarray*} }
\def\eea{\end{eqnarray*} }
\def\be{\begin{equation} }
\def\ee{\end{equation} }

\def\proof{\noindent{\it Proof: }}
\def\qed{\ifhmode\unskip\nobreak\fi\ifmmode\ifinner
\else\hskip5 pt \fi\fi\hbox{\hskip5 pt \vrule width4 pt
height6 pt  depth1.5 pt
\hskip 1pt }}
\begin{document}

\title{ Hypersurfaces with constant sectional curvature
of  $\Sf^n\times \R$ and $\Hy^n\times \R$.}

\author {Fernando Manfio \& Ruy Tojeiro}
\date{}
\maketitle
\addtocounter{equation}{1}

\begin{abstract} We classify the  hypersurfaces of $\Sf^n\times \R$ and $\Hy^n\times \R$ with constant sectional curvature and dimension  $n\geq 3$. \end{abstract}
\medskip
\section[Introduction]{Introduction} The submanifold geometry of the product spaces $\Sf^n\times \R$ and $\Hy^n\times \R$ has been extensively studied in the last years. Here $\Sf^n$ and $\Hy^n$ denote the sphere and hyperbolic space of dimension $n$, respectively. Emphasis has been given on minimal and constant mean curvature surfaces in $\Sf^2\times \R$ and $\Hy^2\times \R$, starting with the work in  \cite{abresch} and \cite{harold}, among others. 
See \cite{mira} for an updated list of references on this topic.

Surfaces of constant {\em Gaussian\/} curvature of $\Sf^2\times \R$ and $\Hy^2\times \R$ were investigated in \cite{aledo} and \cite{espinar}, with special attention to their global properties (see also \cite{irene} for a local study in $\Hy^2\times \R$). In particular, nonexistence of {\em complete\/} surfaces of constant Gaussian curvature $c$  in $\Sf^2\times \R$ (respectively, $\Hy^2\times \R$) was established for $c<-1$ and $0<c<1$ (respectively, $c<-1$). It was also shown that a complete surface of constant Gaussian curvature $c>1$  in $\Sf^2\times \R$ (respectively, $c>0$ in $\Hy^2\times \R$) must be a rotation  surface. Moreover, the profile curves of such surfaces have been explicitly determined.   

Our aim in this paper is to classify all hypersurfaces with constant sectional curvature and dimension $n\geq 3$ of $\Sf^n\times \R$ and $\Hy^n\times \R$. It turns out that for $n\geq 4$ a hypersurface of constant sectional curvature $c$  in $\Sf^n\times \R$ (respectively, $\Hy^n\times \R$) only exists, even locally,  if $c\geq 1$ (respectively, $c\geq -1$), and for any such values of $c$ it  must be an open subset of a  complete rotation hypersurface. In the case $n=3$, exactly one  class of  nonrotational hypersurfaces of  $\Sf^n\times \R$ and  $\Hy^n\times \R$ with constant sectional curvature  arises. Each  hypersurface  in this class in $\Sf^3\times \R$ (respectively, $\Hy^3\times \R$) has constant sectional curvature $c\in (0,1)$ (respectively, $c\in (-1,0)$), and is constructed in an explicit way by means of a family of parallel flat surfaces in  $\Sf^3$ (respectively, $\Hy^3$). An interesting property of such a hypersurface is that its unit normal vector field makes a constant angle with the unit vector field spanning the  factor $\R$.    All  surfaces in $\Sf^2\times \R$ and $\Hy^2\times \R$ with this  property  were classified in \cite{vrancken} and \cite{munteanu}, where they were called {\em constant angle surfaces\/}. Here  we give a simple proof of a generalization of this result to constant angle hypersurfaces of arbitrary dimension of both $\Sf^n\times \R$ and $\Hy^n\times \R$.

\section[Preliminaries]{Preliminaries}

Let $\Q_\e^n$ denote either the sphere $\Sf^n$ or hyperbolic space $\Hy^n$, according as $\e=1$ or $\e=-1$, respectively. In order to study hypersurfaces  $f\colon\,M^n\to \Q_\e^n\times\R$, our approach is to regard $f$ as an isometric immersion into $\Ee^{n+2}$, where $\Ee^{n+2}$ denotes either Euclidean space or Lorentzian space of dimension $(n+2)$, according as $\e=1$ or $\e=-1$, respectively. More precisely, let $(x_1,\ldots,x_{n+2})$ be the standard  coordinates on $\Ee^{n+2}$  with respect to which the flat metric  is written as 
$$d s^2=\e \,d x_1^2+d
x_2^2+\ldots+d x_{n+2}^2.$$
Regard $\Ee^{n+1}$ as 
$$\Ee^{n+1}=\{(x_1,\ldots,x_{n+2})\in\Ee^{n+2}:x_{n+2}=0\}$$
and 
$$ \Q_\e^n=\{(x_1,\ldots,x_{n+1})\in\Ee^{n+1}:\e\,x_1^2+x_2^2+\ldots+x_{n+1}^2=\e\}\,\,\,\mbox{(with $x_1>0$ if $\e=-1$)}.$$
Then we consider the inclusion 
$$
i\colon\, \Q_\e^n\times\R\to\Ee^{n+1}\times\R=\Ee^{n+2}
 $$
 and study the composition $i\circ f$, which we also denote by $f$.

Given a hypersurface $f:M^n\to\Q_\e^n\times\R$, let $N$ denote  a unit normal
vector field to $f$ and let $\frac{\partial}{\partial t}$ be a unit
vector field tangent to the second factor. Then, a vector field  $T$  and a smooth function $\nu$ on $M^n$ are defined by
$$\frac{\partial}{\partial t}=f_* T+\nu N.$$
Notice that $T$ is the gradient of the height function $h=\< f,\frac{\partial}{\partial t}\>$.

Two trivial classes of hypersurfaces of $\Q_\e^n\times \R$ arise if either $\nu$ or $T$  vanishes identically:

\begin{proposition}\label{prop:trivial} Let $f:M^n\to\Q_\e^n\times\R$ be a hypersurface.
\begin{itemize}\item[$(i)$] If $T$  vanishes identically, then $f(M^n)$ is an open subset of a slice $\Q^n_\e\times \{t\}$.
\item[$(ii)$] If $\nu$ vanishes identically, then $f(M^n)$ is an open subset of a Riemannian product $M^{n-1}\times \R$, where $M^{n-1}$ is a hypersurface of $\Q_\e^{n}$.
\end{itemize}
\end{proposition}

Let $\nabla$ and $R$ be the Levi-Civita connection and the curvature
tensor of $M^n$, respectively, and let $A$ be the shape operator of
$f$ with respect to $N$. Then the Gauss and Codazzi equations are 
\be\label{eq:Gauss}\begin{array}{l}
 R(X,Y)Z=(AX\wedge
AY)Z+\e\big((X\wedge Y)Z -\langle
Y,T\rangle(X\wedge T)Z+\langle
X,T\rangle(Y\wedge T)Z\big),
\end{array}\ee
and
\begin{eqnarray}\label{eq:Codazzi}
\nabla_XAY-\nabla_YAX-A[X,Y]=\e\,\nu(X\wedge Y)T,
\end{eqnarray}
respectively, where $X,Y,Z\in TM$. Moreover,  the fact that
$\frac{\partial}{\partial t}$ is parallel in $\Q_\e^n\times\R$ yields for all $X\in TM$ that
\begin{eqnarray}\label{eq:NablaT}
\nabla_XT=\nu AX,
\end{eqnarray}
and 
\begin{eqnarray}\label{eq:DerivadaNu}
X(\nu)=-\langle AX,T\rangle.
\end{eqnarray}

\section[A basic lemma]{A basic lemma}

Our main goal in this section is to prove the following lemma.

\begin{lemma}\label{le:Tprincipal1}
Let $f:M_c^n\to\Q_\e^n\times\R$ be a hypersurface of dimension $n\geq 3$ and constant sectional curvature $c\neq 0$. Assume that $T\neq 0$ at $x\in M^n_c$. Then $T$ is a principal direction at $x$. 
\end{lemma}

Lemma \ref{le:Tprincipal1} will follow by putting together Lemma \ref{lem:Rperpzero} and Proposition \ref{prop:fnb} below:

\begin{lemma}\label{lem:Rperpzero}  Let $f:M^n\to\Q_\e^n\times\R$ be a hypersurface. Suppose that $T\neq 0$  at $x\in M^n$. Then $f$ has flat normal bundle at $x$ as an isometric immersion into $\Ee^{n+2}$ if and only if $T$ is a principal direction at $x$.
\end{lemma}

\begin{proposition}\label{prop:fnb}
Any isometric immersion $g:M_c^n\to\Ee^{n+2}$ of a Riemannian manifold with dimension $n\geq 3$ and constant sectional curvature $c\neq 0$ has flat normal bundle.  
\end{proposition}

Lemma \ref{lem:Rperpzero} was first proved in \cite{dillen} for $n=2$ and $\e=1$. A proof of the general case can be found in \cite{tojeiro}. For the proof of Proposition \ref{prop:fnb} we make use of standard facts from \cite{moore} on the theory of flat bilinear forms. Recall that a symmetric bilinear form $\beta\colon\, V\times V\to W$, where $V$ and $W$ are finite-dimensional vector spaces,  is said to be \emph{flat} with respect to an inner product $\<\,\,,\,\,\>\colon\,W\times W\to \R$ if  
$$
\<\beta(X,Y),\beta(Z,T)\>-\<\beta(X,T),\beta(Z,Y)\>=0
$$
for all  $X,Y,Z,T\in V$. Clearly, the standard example of a flat bilinear form is the second fundamental form of an isometric immersion between space forms  with the same constant sectional curvature. 

Denote by $N(\beta)\subset V$ the \emph{nullity subspace } of $\beta$, given by
$$
N(\beta)=\{X\in V\colon \beta(X,Y)=0 : Y\in V\},
$$
and by $S(\beta)\subset W$ its \emph{image subspace}
$$
S(\beta)=\spa\{\beta(X,Y) : X, Y\in V\}.
$$
 The next result is a basic fact on flat bilinear forms (cf. Corollary $1$ and Corollary $2$ in \cite{moore}) : 

\begin{theorem} \label{thm:moore1} {\em \cite{moore}} Let $\beta\colon V\times V\to W$ be a flat bilinear form with respect to  an inner product $\<\,\,,\,\,\>$ on $W$. Assume that $\<\,\,,\,\,\>$ is either positive-definite or Lorentzian and, in the latter case, suppose that $S(\beta)$ is a nongenerate subspace of $W$, i.e., $S(\beta)\cap S(\beta)^\perp=\{0\}$. Then  
$$
\dim N(\beta)\ge\dim V-\dim S(\beta).
$$
\end{theorem}

Another fact we will need  in order to handle the case $n=3$ in Proposition \ref{prop:fnb} is the following  consequence of Theorem $2$ in  \cite{moore}:

\begin{theorem} \label{thm:moore2} {\em \cite{moore}} Let $\beta\colon V\times V\to W$ be a flat bilinear form with respect to  an inner product $\<\,\,,\,\,\>$ on $W$. Assume that $\dim V=\dim W$, that $N(\beta)=\{0\}$ and that $\<\,\,,\,\,\>$ is either positive-definite or Lorentzian. Moreover, in the latter case suppose that there exists a vector $e\in W$ such that $\<\beta(\,\,,\,\,),e\>$ is positive definite. Then there exists a diagonalyzing basis  $\{e_1,\ldots, e_n\}$ for $\beta$, i.e., $\beta(e_i,e_j)=0$ for $1\leq i\neq j\leq n$.
\end{theorem}

\noindent {\em Proof of Proposition  \ref{prop:fnb}:}
First recall that $\R^{n+2}$ admits an umbilical inclusion  $i$ into both hyperbolic space $\Hy_c^{n+3}$ and the 
Lorentzian sphere $\Sf_c^{n+2,1}$ of constant sectional curvature $c$, according as $c<0$ or $c>0$, respectively, i.e., its second fundamental form $\alpha$ is
\[\alpha(X,Y)=\sqrt{|c|}\langle X,Y\rangle\eta,\]
where $\eta$ is one of the two normal vectors such that
$\langle\eta,\eta\rangle=-\mbox{sgn(c)}$, where $\mbox{sgn(c)}=c/|c|$. Similarly, Lorentzian space $\Le^{n+2}$ admits  umbilical inclusions into $\Hy_c^{n+2,1}$ or  $\Sf_c^{n+1,2}$, according as $c<0$ or $c>0$, respectively. 

Then, the second fundamental form $\alpha_\phi=g^*\alpha+i_*\alpha_g$ of
$\phi=i\circ g$ at every $x\in M_c^n$ is a flat bilinear form with respect to the inner product $\<\,\,,\,\,\>$ on its three-dimensional normal space. The inner product  $\<\,\,,\,\,\>$ is positive-definite if $c<0$ and $\Ee^{n+2}=\R^{n+2}$, Lorentzian if either $c>0$ and $\Ee^{n+2}=\R^{n+2}$ or if $c<0$ and $\Ee^{n+2}=\Le^{n+2}$, and  has index two if $c>0$ and $\Ee^{n+2}=\Le^{n+2}$. In the latter case, $\alpha_\phi$ is also flat with respect to the Lorentzian inner product  $-\<\,\,,\,\,\>$. Moreover, since 
$$\langle\alpha_\phi(,),i_*\eta\rangle=\langle\alpha(,),\eta\rangle=-\mbox{sgn(c)}\sqrt{|c|}\<\,\,,\,\,\>,$$  it follows that  $N(\alpha_\phi)=\{0\}$.  Let us consider the two possible cases:\vspace{2ex}\\
\noindent {\em $(i)$ $S(\alpha_\phi)$ is nondegenerate :\/}
in this case Theorem \ref{thm:moore1} gives
\[\dim S(\alpha_\phi)\geq n-\dim N(\alpha_\phi)=n.\]
Since $\dim S(\alpha_\phi)\leq 3$, this implies that $n=3=\dim S(\alpha_\phi)$. Since $\langle\alpha_\phi(,),-\mbox{sgn(c)}i_*\eta\rangle$ is positive definite, it
follows from Theorem \ref{thm:moore2}  that there exists a basis $\{e_1,\ldots, e_n\}$ of $T_x M^n_c$ such that
$\alpha_\phi(e_i,e_j)=0$ for $i\neq j$. In particular, we have 
$$0=\<\alpha_\phi(e_i,e_j),i_*\eta\>=-\mbox{sgn(c)}\sqrt{|c|}\langle e_i,e_j\rangle\,\,\,\mbox{for}\,\,\,i\neq j,$$ that is, $\{e_1,\ldots, e_n\}$ is an orthogonal basis. Since $\{e_1,\ldots, e_n\}$ also diagonalizes $\alpha_{g}$, we conclude that $g$ has flat normal bundle. 
 \vspace{1ex}\\
\noindent {\em $(ii)$ $S(\alpha_\phi)$ is degenerate :\/}
in this case, there exists a nonzero vector $\rho\in
S(\alpha_\phi)\cap S(\alpha_\phi)^\perp$. Writing
$\rho=\eta+i_*\zeta$, with $\zeta$ a unit normal vector to $g$, we obtain from $0=\langle\alpha_\phi(X,Y),\rho\rangle$ for all $X,Y\in T_xM_c^n$   that
\[\langle\alpha_{g}(X,Y),\zeta\rangle=\mbox{sgn(c)}\sqrt{|c|}\langle X,Y\rangle,\]
for all $X,Y\in T_xM_c^n$, i.e., $g$ has an  umbilical
normal direction. Since $g$ has codimension two, the Ricci equation  implies that its normal bundle is flat. \qed\vspace{.2cm}\\

The flat case $c=0$ can also be handled by means of Theorem \ref{thm:moore1}:

\begin{lemma}\label{le:Tprincipal2}
Let $f:M_0^n\to\Q_\e^n\times\R$ be a flat hypersurface of dimension $n\geq 3$.  Assume that $T\neq 0$ at $x\in M^n_0$.
\begin{itemize} \item[$(i)$] 
If $\e=1$, then $\nu$ vanishes at $x$.
\item[$(ii)$] If $\e=-1$, then  either $\nu$ vanishes at $x$ or $A_N=A_\xi$ for one of the two possible choices of a unit normal vector $N$ to $f$ at $x$. 
\end{itemize}
In any case, $T$ is a principal direction of $f$ at $x$.
\end{lemma}
\proof Regard $f$ as an isometric immersion into $\Ee^{n+2}$. 
Then, its  second fundamental form $\a$  is a flat bilinear map by the Gauss equation.  Let $\xi$ denote  the outward pointing unit normal vector field to $\Q_\e^n\times\R$. Then it is easily seen that the shape operator of $ f$ with respect to $\xi$ is given by
\be\label{eq:eigenvectorAxi}
A_\xi T=-\nu^2T\,\,\,\,\mbox{and}\,\,\,\,A_\xi X=-X\,\,\,\mbox{for}\,\,X\in \{T\}^\perp.
\ee
Assume that  $\nu\neq 0$ at $x\in M_0^n$. Then $A_\xi$, and hence  $\alpha$, has trivial kernel by (\ref{eq:eigenvectorAxi}).   If $\e=1$, it follows from Theorem \ref{thm:moore1} that 
$$2\geq \dim S(\alpha)\geq n,$$ a contradiction that proves $(i)$.
If $\e=-1$, Theorem \ref{thm:moore1} in the Lorentzian case implies that $S(\alpha)$ is a degenerate subspace of the two-dimensional normal space of $f$ in $\Ee^{n+2}$ at $x$. Hence $S(\alpha)$ is spanned by the light-like vector $i_*N+\xi$ for one of the two unit normal vectors $N$ to $f$ in $\Q_\e^n\times \R$  at $x$. But the fact that $i_*N+\xi \in S(\alpha)^\perp$ just means that $A_N=A_\xi$. 

For the last assertion, notice that a point where $\nu$ vanishes is a local minimum for $\nu$, hence  $A_NT=0$ at $x$ by (\ref{eq:DerivadaNu}).
\qed

\section[Rotation hypersurfaces]{Rotation hypersurfaces}

Rotation hypersurfaces of $\Sf^n\times\R$ and $\Hy^n\times\R$ have been defined and their principal curvatures computed in \cite{veken}, as an extension of the work in \cite{dcd} on rotation hypersurfaces of space forms.

With notations as in Section $2$, let $P^3$ be a three-dimensional subspace of $\Ee^{n+2}$ containing the $\frac{\partial}{\partial x_1}$ and the 
$\frac{\partial}{\partial x_{n+2}}$ directions. Then $(\Q_\e^n\times\R)\cap P^3=
\Q_\e^1\times\R$.  Denote by $\mathcal I$
the group of isometries of $\Ee^{n+2}$ that fix pointwise a two-dimensional subspace $P^2\subset P^3$ also
containing the $\frac{\partial}{\partial x_{n+2}}$-direction. 
Consider a curve  $\alpha$  in
$\Q_\e^1\times\R\subset P^3$ that lies in one of the two half-spaces of $P^3$ determined by  $P^2$.

\begin{definition} {\em A {\em rotation
hypersurface  in $\Q_\e^n\times\R$ with profile curve $\alpha$ and axis $P^2$}
is the orbit of $\alpha$ under the action of $\mathcal I$.}
\end{definition}

We will always assume that $P^3$ is spanned by $\frac{\partial}{\partial x_1}$, $\frac{\partial}{\partial x_{n+1}}$  and 
$\frac{\partial}{\partial x_{n+2}}$. In the case $\e=1$, we also assume that $P^2$ is spanned by $\frac{\partial}{\partial x_1}$ and 
$\frac{\partial}{\partial x_{n+2}}$, and that the curve $\alpha$ is parametrized by arc length as
$$\alpha(s)=(\sin(k(s)),0,\ldots,0,\cos(k(s)),h(s)),$$
where  $s$ runs over an interval $I$ where $\cos (k(s))\geq 0$, so that $\alpha(I)$ is contained in  a closed half-space determined by  $P^2$. Here $k,h\colon\,I\to \R$ are smooth functions satisfying
\be\label{eq:arclength}k'(s)^2+h'(s)^2=1\,\,\,\mbox {for all $s\in I$}.\ee
 In this case, the  rotation
hypersurface  in $\Sf^n\times\R$ with profile curve $\alpha$ and axis $P^2$ can be parametrized by
\be\label{eq:parametrization} f(s,t)=(\sin(k(s)),\cos(k(s))\varphi_1(t),\ldots,\cos(k(s))\varphi_n(t),h(s)),\ee
where $t=(t_1,\ldots,t_{n-1})$ and
$\varphi=(\varphi_1,\ldots,\varphi_n)$  parametrizes
 $\Sf^{n-1}\subset\R^n$. The metric induced by $f$ is
\be\label{eq:met1} d\sigma^2=ds^2+\cos^2(k(s))dt^2,\ee
where $dt^2$ is the standard metric of $\Sf^{n-1}$. 
 
 For $\e=-1$, one has three distinct possibilities, according as $P^2$ is Lorentzian, Riemannian or degenerate, respectively. We call $f$, accordingly, a rotation hypersurface of {\em spherical\/}, {\em hyperbolic\/} or {\em parabolic\/} type, because the orbits of $\mathcal I$ are spheres, hyperbolic spaces or horospheres, respectively.  In the first case, we can assume that 
 $P^2$ is spanned by $\frac{\partial}{\partial x_1}$ and 
$\frac{\partial}{\partial x_{n+2}}$ and that the curve $\alpha$ is parametrized by 
\be\label{eq:paramcurve}\alpha(s)=(\cosh(k(s)),0,\ldots,0,\sinh(k(s)),h(s)).\ee
Then $f$ can be parametrized by
\be\label{eq:parametrization2} f(s,t)=(\cosh(k(s)),\sinh(k(s))\varphi_1(t),\ldots,\sinh(k(s))\varphi_n(t),h(s)).\ee
The induced metric is 
 \be\label{eq:met2} d\sigma^2=ds^2+\sinh^2(k(s))dt^2,\ee
where $dt^2$ is the standard metric of $\Sf^{n-1}$.

In the second case, assuming  that 
 $P^2$ is spanned by $\frac{\partial}{\partial x_{n+1}}$ and 
$\frac{\partial}{\partial x_{n+2}}$,  the curve $\alpha$ can also be parametrized as in  (\ref{eq:paramcurve}), and a  parametrization of $f$ is
\be\label{eq:parametrization3} f(s,t)=(\cosh(k(s))\varphi_1(t),\ldots,\cosh(k(s))\varphi_n(t),\sinh( k(s)),h(s)),\ee
where $t=(t_1,\ldots,t_{n-1})$ and
$\varphi=(\varphi_1,\ldots,\varphi_n)$  parametrizes
 $\Hy^{n-1}\subset\Le^n$. The induced metric is 
 \be\label{eq:met3} d\sigma^2=ds^2+\cosh^2(k(s))dt^2,\ee
where $dt^2$ is the standard metric of $\Hy^{n-1}$.

 Finally, when $P^2$ is degenerate, we choose a pseudo-orthonormal basis 
 $$e_1=\frac{1}{\sqrt{2}}\left(-\frac{\partial}{\partial x_{1}}+\frac{\partial}{\partial x_{n+1}}\right),\,\,\,\, e_{n+1}=\frac{1}{\sqrt{2}}\left(\frac{\partial}{\partial x_{1}}+\frac{\partial}{\partial x_{n+1}}\right), \,\,\,e_j=\frac{\partial}{\partial x_{j}}, $$
 for $j\in \{2,\ldots, n,n+2\}$, and assume that  $P^2$ is spanned by $e_{n+1}$ and $e_{n+2}$. Notice that $\<e_1,e_1\>=0=\<e_{n+1},e_{n+1}\>$ and $\<e_1,e_{n+1}\>=1$. Then, we can parametrize $\alpha$ by 
 $$\alpha(s)=\left(k(s),0,\ldots,0,-\frac{1}{2k(s)},h(s)\right),$$
 with \be\label{eq:arclength2} k(s)> 0\,\,\,\,\mbox{and}\,\,\,\,(\ln k)'^{\,2}(s)+h'(s)^2=1,\ee
 and a parametrization of $f$ is 
 \be\label{eq:parametrization4} f(s,t_2,\ldots, t_n)=\left(k(s),k(s)t_2,\ldots, k(s)t_n, -\frac{1}{2k(s)}-\frac{k(s)}{2}\sum_{i=2}^nt_i^2, h(s)\right),\ee
 whose induced metric is
 \be\label{eq:met4} d\sigma^2=ds^2+k^2(s)dt^2,\ee
where $dt^2$ is the standard metric of $\R^{n-1}$. 

\begin{remark}{\em Our definition of a rotation hypersurface in $\Q^n_\e\times \R$ was taken from \cite{veken}, and it naturally extends the one given in \cite{dcd} for space forms. For $\e=-1$, it differs from that used in \cite{aledo}, where only rotation surfaces of spherical type were considered.}
\end{remark} 

We are now in a position to classify rotation hypersurfaces of $\Q^n_\e\times \R$ with constant sectional curvature $c$ and dimension $n\geq 3$. We state separately the cases $\e=1$ and $\e=-1$:

\begin{theorem}\po\label{thm:rotationcsc1} Let  $f\colon\, M_c^n\to \Sf^n\times \R$  be a rotation hypersurface with constant sectional curvature $c$ and dimension $n\geq 3$. Then $c\geq 1$. Moreover,  
 \begin{itemize}
 \item[$(i)$] if $c=1$ then $f(M_c^n)$ is an open subset of a slice $\Sf^{n}\times \{t\}$.
 \item[$(ii)$] if  $c> 1$ then $f(M_c^n)$ is an open subset of a complete hypersurface that can be parametrized by (\ref{eq:parametrization}), with
 \be\label{eq:k1}
k(s)=\arccos\left(\frac{1}{\sqrt{c}}\sin (\sqrt{c}\,s)\right)
\ee
and
\be\label{eq:h1}
h(s)=-{\sqrt\frac{c-1}{c}}\ln\left(\frac{\cos(\sqrt{c}\,s)+\sqrt{c-\sin^2(\sqrt{c}\, s)}}{1+\sqrt{c}}\right),\,\,\,\,s\in [0,\pi/\sqrt{c}].
\ee
 \end{itemize}
 \end{theorem}

 \begin{theorem}\po\label{thm:rotationcsc2} Let  $f\colon\, M_c^n\to \Hy^n\times \R$  be a rotation hypersurface with constant sectional curvature $c$ and dimension $n\geq 3$. Then $c\geq -1$. Moreover, 
 \begin{itemize}
 \item[$(i)$] if $c=-1$ then $f(M^n)$ is an open subset of a slice $\Hy^{n}\times \{t\}$.
 \item[$(ii)$] if $c\in (-1,0)$ then  one of the following possibilities holds:
 \begin{itemize}
 \item[$(a)$] $f(M^n)$ is an open subset of a complete hypersurface of spherical type that can be parametrized by (\ref{eq:parametrization2}), with
 \be\label{eq:k2}
k(s)=\mbox{\em arcsinh}\left(\frac{1}{\sqrt{-c}}\sinh (\sqrt{-c}\,s)\right)
\ee
and
\be\label{eq:h2}
h(s)=\sqrt{\frac{c+1}{-c}}\ln\left(\frac{\cosh(\sqrt{-c}\,s)+\sqrt{-c+\sinh^2(\sqrt{-c}\, s)}}{1+\sqrt{-c}}\right).
\ee
 \item[$(b)$] $f(M^n)$ is an open subset of a  complete hypersurface of hyperbolical type that can be parametrized by (\ref{eq:parametrization3}), with
\be\label{eq:k3}
 k(s)=\mbox{\em arccosh} \frac{1}{\sqrt{-c}}\cosh(\sqrt{-c}\,s)
 \ee
 and 
 \be\label{eq:h3}
 h(s)=\sqrt{\frac{c+1}{-c}}\ln\left(\sinh(\sqrt{-c}\,s)+\sqrt{c+\cosh^2(\sqrt{-c}\, s)}\right).
\ee
 \item[$(c)$] $f(M^n)$ is an open subset of a complete hypersurface of parabolical type that can be parametrized by (\ref{eq:parametrization4}), with
 \be\label{eq:k4}
 k(s)=\exp  \sqrt{-c}\, s
 \ee
 and 
 \be\label{eq:h4}
h(s)=\sqrt{1+c}\,s.
\ee
\end{itemize}
\item[$(iii)$] if $c=0$, then one of the following possibilities holds:
 \begin{itemize}
 \item[$(a)$] $f(M^n)$ is an open subset of a complete hypersurface of spherical type that can be parametrized by (\ref{eq:parametrization2}), with
 \be\label{eq:k5}
k(s)=\mbox{\em arcsinh} (s)
\ee
and
\be\label{eq:h5}
h(s)=-1+\sqrt{1+s^2}.
\ee
 \item[$(b)$] $f(M^n)$ is an open subset of a Riemannian product $M^{n-1}\times \R$, where $M^{n-1}$ is a horosphere of $\Hy^n$.
 \end{itemize}
 \item[$(iv)$] if $c>0$,  then $f(M^n)$ is an open subset of a complete hypersurface of spherical type that can be parametrized by (\ref{eq:parametrization2}), with
 \be\label{eq:k6}
k(s)=\mbox{\em arcsinh}\left(\frac{1}{\sqrt{c}}\sin (\sqrt{c}\,s)\right)
\ee
and
\be\label{eq:h6}
h(s)=-{\sqrt\frac{c+1}{c}}\arctan\left(\frac{\cos(\sqrt{c}\,s)}{\sqrt{c+\sin^2(\sqrt{c}\, s)}}\right).
\ee
 \end{itemize}
 \end{theorem}
 
 \begin{remark}{\em The hypersurfaces in Theorems \ref{thm:rotationcsc1} and \ref{thm:rotationcsc2} also occur in dimension $n=2$. In particular, those in parts $(ii)-b)$ and $(ii)-c)$ of Theorem \ref{thm:rotationcsc2}  provide examples of complete surfaces of constant Gaussian curvature $c\in (-1,0)$ in $\Hy^2\times \R$ that do not appear in \cite{aledo}.}
 \end{remark}  
 
 For the proof of Theorems \ref{thm:rotationcsc1} and \ref{thm:rotationcsc2} we make use of the following  fact:

 \begin{proposition}\po\label{prop:warped} Assume that the  warped product $I\times_\rho \Q_\delta^n$, $n\geq 2$, $\delta\in \{-1,0,1\}$, has constant sectional curvature $c$. 
 \begin{itemize}
 \item[$(i)$] If $c>0$, then $\delta=1$ and
 $\rho(s)=\frac{1}{\sqrt{c}}\sin(\sqrt{c}\,s+\theta_0),\,\,\,\theta_0\in \R.$
 \item[$(ii)$] If $c=0$, then one of the following possibilities holds: 
  \begin{itemize}
  \item[$(a)$] $\delta=1$ and $\rho(s)=\pm s+s_0,\,\,\,s_0\in \R.$
  \item[$(b)$] $\delta=0$ and $\rho(s)=A\in \R$.
   \end{itemize}
  \item[$(iii)$]  If $c<0$, then one of the following possibilities holds: 
  \begin{itemize}
  \item[$(a)$] $\delta=-1$ and $\rho(s)=\frac{1}{\sqrt{-c}}\cosh(\sqrt{-c}\,s+\theta_0),\,\,\,\theta_0\in \R.$
  \item[$(b)$] $\delta=0$ and $\rho(s)=\exp(\pm \sqrt{-c}\, s+s_0),\,\,\,s_0\in \R$.
  \item[$(c)$] $\delta=1$ and $\rho(s)=\frac{1}{\sqrt{-c}}\sinh(\sqrt{-c}\,s+\theta_0),\,\,\,\theta_0\in \R.$
 \end{itemize}
 \end{itemize}
 \end{proposition}
 \proof  In a warped product $I\times_\rho \Q_\delta^n$, $n\geq 2$, the sectional curvature along a plane tangent to $\Q_\delta^n$ is $(\delta-(\rho')^2)/\rho^2$, whereas the sectional curvature along a plane  spanned by unit vectors $\d/\d s$ and $X$ tangent to $I$ and $\Q_\delta^n$, respectively, is $-\rho''/\rho$. Therefore, $I\times_\rho \Q_\delta^n$ has constant sectional curvature $c$ if and only if 
 \be\label{eq:rho} (\rho')^2+c\rho^2=\delta.\ee
  Notice that $-\rho''/\rho=c$, or equivalently, 
  \be\label{eq:rho2}\rho''+c\rho=0,\ee follows by differentiating (\ref{eq:rho}).  If $c>0$, we obtain from (\ref{eq:rho}) that $\delta=1$. Moreover, by (\ref{eq:rho2}) we have that
  $$\rho(s)=A\cos \sqrt{c}\,s+B\sin \sqrt{c}\,s$$
  for some $A,B\in \R$, which gives $(\rho')^2+c\rho^2=c(A^2+B^2)$.  From (\ref{eq:rho}) we get $c(A^2+B^2)=1$, hence we may write $$A=\frac{1}{\sqrt{c}}\sin \theta_0\,\,\,\,\mbox{and}\,\,\,\,B=\frac{1}{\sqrt{c}}\cos \theta_0$$
  for some $\theta_0\in \R$. It follows that
  $$\rho(s)=\frac{1}{\sqrt{c}}\sin (\sqrt{c}\,s+\theta_0).$$
  The remaining cases are similar.\vspace{2ex}\qed

 \noindent {\em Proof of Theorems \ref{thm:rotationcsc1} and \ref{thm:rotationcsc2}:} First we determine the possible values of $c$ for a  rotation hypersurface $f\colon\, M_c^n\to \Q^n_\e\times \R$  with constant sectional curvature $c$ and dimension $n\geq 3$. If $T$ vanishes on an open subset, then $c=\e$ by Proposition \ref{prop:trivial}. Otherwise, we can assume that $T$ is nowhere vanishing. Then $f$ has exactly two distinct principal curvatures $\lambda$ and $\mu\neq 0$, the first one being simple with $T$ as principal direction (cf. \cite{veken}). Let $\{T,X_1,\ldots,X_{n-1}\}$ be an orthogonal basis of
eigenvectors of $A$ at $x$, with $$AT=\lambda T\,\,\,\mbox{ and}\,\,\,AX_i=\mu X_i,\,\,\,1\leq i\leq n-1.$$   From the Gauss equation (\ref{eq:Gauss}) of $f$ for $X=X_i$ and $Y=Z=X_j$, $i\neq j$, we get
$$
c-\e=\mu^2,
$$
 and hence $c>\e$. This proves the first assertions in Theorems \ref{thm:rotationcsc1} and \ref{thm:rotationcsc2}.

 Now assume that $\e=1$. Then $f$ can be parametrized by (\ref{eq:parametrization}), with $k(s)$ and $h(s)$ satisfying (\ref{eq:arclength}), and  the metric induced by $f$ is given by (\ref{eq:met1}). Since $c\geq 1$, by Proposition~\ref{prop:warped} we must have 
 $$\cos(k(s))=\frac{1}{\sqrt{c}}\sin(\sqrt{c}\,s+\theta_0)$$
 for some $\theta_0\in \R.$  Replacing $s$ by $s-\theta_0/\sqrt{c}$, we can assume that $\theta_0=0$. 
If $c=1$, then $f$  just parametrizes an open subset of a slice $\Sf^{n}\times \{t\}$. If $c> 1$, we obtain that $k(s)$ and $h(s)$ are given by (\ref{eq:k1}) and (\ref{eq:h1}), respectively. The corresponding profile curve is exactly that of the complete surface of constant sectional curvature $c$ in $\Sf^2\times \R$ determined in \cite{aledo}, and their argument also applies to show the completeness of $f$ in any dimension $n\geq 3$.\vspace{1ex}

From now on we deal with the case $\e=-1$. Assume first that $f$ is of spherical type. Then $f$ can be parametrized by 
(\ref{eq:parametrization2}), with $k(s)$ and $h(s)$ satisfying (\ref{eq:arclength}), and  the metric induced by $f$ is given by (\ref{eq:met2}). By Proposition \ref{prop:warped}, the warping function $\sinh(k(s))$ must be equal to 
$$ \frac{1}{\sqrt{c}}\sin(\sqrt{c}\,s+\theta_0),\,\,\,\,\frac{1}{\sqrt{-c}}\sinh(\sqrt{-c}\,s+\theta_0),\,\,\,\theta_0\in \R,\,\,\,\,\mbox{or}\,\,\,\,\,\pm s+s_0,\,\,\,s_0\in \R,$$
according as $c>0$, $c<0$ or $c=0$, respectively. After suitably replacing the parameter $s$, we can assume that $\theta_0=0$ in the first two cases, and that $\sinh(k(s))=s$ in the last one. Each possibility gives rise to the expressions (\ref{eq:k2}), (\ref{eq:k6}) and (\ref{eq:k5}) for $k(s)$, and (\ref{eq:h2}), (\ref{eq:h6}) and (\ref{eq:h5}) for $h(s)$, respectively. The corresponding profile curves are exactly those of the complete rotation surfaces with constant sectional curvature of spherical type determined in \cite{aledo}, and the completeness of the corresponding hypersurfaces can be seen in the same way as in \cite{aledo}.

Now suppose that $f$ is of hyperbolical type. Then, it can be parametrized by 
(\ref{eq:parametrization3}), with $k(s)$ and $h(s)$ satisfying (\ref{eq:arclength}), and  the  induced metric is (\ref{eq:met3}). Since $c\geq -1$, by Proposition \ref{prop:warped} we must have $c\in [-1, 0)$ and 
$$\cosh(k(s))=\frac{1}{\sqrt{-c}}\cosh(\sqrt{-c}\,s+\theta_0),\,\,\,\theta_0\in \R.$$
As before, we can assume that $\theta_0=0$. If $c=-1$, then $f(M^n)$ is an open subset of a slice $\Hy^{n}\times \{t\}$. Otherwise,  $k$ and $h$ are  given by (\ref{eq:k3}) and (\ref{eq:h3}), respectively.

Finally, suppose that   $f$ is of parabolical type. Then, it can be parametrized by 
(\ref{eq:parametrization4}), with $k(s)$ and $h(s)$ satisfying (\ref{eq:arclength2}), and  the  induced metric is (\ref{eq:met4}). By Proposition \ref{prop:warped}, we must have $c\leq 0$ and 
$$k(s)=A\in \R\,\,\,\,\,\mbox{or}\,\,\,\,\,k(s)=\exp(\pm \sqrt{-c}\, s+s_0),\,\,\,s_0\in \R,$$
according as $c=0$ or $c<0$, respectively. In the first case, $f$ just parametrizes an open subset of a Riemannian product $M^{n-1}\times \R$, where $M^{n-1}$ is a horosphere of $\Hy^n$. In the second case, we can assume that $k(s)=\exp \sqrt{-c}\, s$ and then  $h$ is given by (\ref{eq:h4}). Completeness of the hypersurfaces  in this and the preceding case is straightforward. \vspace{2ex}\qed

\section[Constant angle  hypersurfaces]{Constant angle hypersurfaces}

Let $g\colon\, M^{n-1}\to \Q_\e^n$ be a hypersurface and let $g_s \colon\, M^{n-1}\to \Q_\e^n$ be the family of parallel hypersurfaces to $g$, that is, 
\be\label{eq:parallel}g_s(x)=C_\e(s)g(x)+S_\e(s)N(x),\ee
where $N$ is a unit normal vector field to $g$, 
$$
S_\e(s)=\left\{\begin{array}{l}
\cos s, \,\,\,\mbox{if}\,\,\e=1
\vspace{1.5ex}\\
\cosh s, \,\,\,\mbox{if}\,\,\e=-1
\end{array}\right.\,\,\,\,\,\,\,\,\mbox{and}\,\,\,\,\,\,\,\,\,\,\,
S_\e(s)=\left\{\begin{array}{l}
\sin s, \,\,\,\mbox{if}\,\,\e=1
\vspace{1.5ex}\\
\sinh s, \,\,\,\mbox{if}\,\,\e=-1.
\end{array}\right.
$$
For $\e=1$, write the principal curvatures of $g$ 
as 
$$\lambda_i=\cot \theta_i,\,\,\,\,0<\theta_i<\pi,\,\,\,\,1\leq i\leq m,$$
where the $\theta_i$ form an increasing sequence. For $X$ in the  eigenspace of the shape operator $A_N$ of $g$ corresponding to the principal curvature $\lambda_i$, $1\leq i\leq m$, we have
$${g_s}_*X=g_*(\cos s \,X-\sin s\,A_N X)=(\cos s -\sin s \cot\theta_i)X=\frac{\sin(\theta_i-s)}{\sin \theta_i}X,$$
 Thus, $g_s$ is an immersion at $x$ if and only if $s\neq \theta_i(x)\mbox{(mod $\pi$)}$ for any $1\leq i \leq m$. 

For $\e=-1$, write the principal curvatures of $g$ with absolute value greater than $1$ as 
$$\lambda_i=\coth \theta_i,\,\,\,\,\theta_i\neq 0,\,\,\,\,1\leq i\leq m.$$
 As in the preceding case, for $X$ in the  eigenspace of the shape operator $A_N$ corresponding to the principal curvature $\lambda_i$, $1\leq i\leq m$, we have
$${g_s}_*X=\frac{\sinh(\theta_i-s)}{\sinh \theta_i}X,$$
 Thus, $g_s$ is an immersion at $x$ if and only if $s\neq \theta_i(x)$ for any $1\leq i \leq m$. 

In the case $\e=1$, set
\be\label{eq:U} U:=\{(x,s)\in M^{n-1}\times \R\,:\,s\in (\theta_m(x)-\pi, \theta_1(x))\}.\ee
For $\e=-1$,  let $\theta_+$ (respectively, $\theta_-$) be the least (respectively, greater) of the $\theta_i$ that is  greater than $1$ (respectively, less than $-1$), and set 
\be\label{eq:U} U:=\{(x,s)\in M^{n-1}\times \R\,:\,s\in (\theta_-(x), \theta_+(x))\}.\ee
In both cases, if  $V\subset M^{n-1}$ is an open subset and $I$ is an open interval containing $0$ such that $V\times I\subset U$, then $g_s$ is an immersion on $V$ for every $s\in I$, with 
\be\label{eq:normalparallel}N_s(x)=-\e S_\e(s)g(x)+C_\e(s)N(x)\ee
as a unit normal vector at $x$.

 Now define
$$f\colon\, M^n:=V\times I\to \Q_e^n\times \R\subset \Ee^{n+2}$$  by
\be\label{eq:constantangle}f(x,s)=g_s(x)+ Bs\frac{\d}{\d t},\,\,\,\,B> 0.\ee
Then 
$$f_*X={g_s}_*X,\,\,\,\,\,\mbox{for any}\,\, X\in TM^{n-1},$$
and 
$$f_*\frac{\d}{\d s}=N_s+B\frac{\d}{\d t},$$
where \be\label{eq:normalparallel}N_s(x)=-\e S_\e(s)g(x)+C_\e(s)N(x).\ee
Since $g_s$ is an immersion on $V$ for every $s\in I$, it follows that $f$ is  an immersion on $M^n$ with 
\be\label{eq:eta}\eta(x,s)=-\frac{B}{a}N_s(x)+\frac{1}{a}\frac{\d}{\d t},\,\,\,a=\sqrt{1+B^2}\ee
as a unit normal vector field. Thus, $f$ has the property that
 $$\<\eta,\frac{\d}{\d t}\>=\frac{1}{a}$$
is constant on $M^n$. Following \cite{vrancken},  $f$ was called in \cite{tojeiro} a {\em constant angle hypersurface}. Constant angle surfaces in $\Sf^2\times \R$ and $\Hy^2\times \R$ have been classified in \cite{vrancken} and \cite{munteanu}, respectively. The next result was obtained in \cite{tojeiro} as a consequence of a more general theorem. For the sake of completeness we provide here a simple and direct proof. 

 \begin{theorem}\po\label{thm:constantangle} Any constant angle hypersurface $f\colon\, M^n\to \Q_\e^n\times \R$  is  either an open subset of a  slice $\Q_\e^n\times \{t_0\}$ for some $t_0\in \R$, an open subset of a product $M^{n-1}\times \R$, where $M^{n-1}$ is a hypersurface of $\Q_\e^n$, or it is locally given by the preceding construction.
\end{theorem}
\proof Let $\eta$ be a unit normal vector field to $f$. By assumption, $\nu=\<\eta,\d/\d t\>$ is  a constant  on $M^n$, which we can assume to belong to $[0,1]$. Since $\|T\|^2+\nu^2=1$, the vector field $T$ has also constant length. 
By Proposition \ref{prop:trivial}, the cases $\nu=1$ and $\nu=0$ correspond to the first two possibilities in the statement, respectively. From now on, we assume that $\nu\in (0,1)$, hence $T$ is a vector field whose length is also a constant in $(0,1)$. Since $T$ is a gradient  vector field, its integral curves  are (not unit-speed) geodesics in $M^n$. The fact that $T$  is a gradient also implies that the orthogonal distribution $\{T\}^\perp$ is integrable. Thus, there exists locally a diffeomorphism $\psi\colon\,M^{n-1}\times I\to M^n$, where $I$ is an open interval containing $0$, such that $\psi(x,\cdot)\colon\, I\to M^n$ are  integral curves of $T$ and $\psi(\cdot, s)\colon\, M^{n-1}\to M^n$ are integral manifolds of $\{T\}^\perp$. Set $F= f\circ \psi$, with $f$ being regarded as an isometric immersion into $\Ee^{n+2}$. Then $$X\<F, \frac{\d}{\d t}\>=\< f_* \psi_*X, \frac{\d}{\d t}\>=\<\psi_*X, T\>=0$$
for any $X\in TM^{n-1}$. Thus $\<F(x,s), \frac{\d}{\d t}\>=\rho(s)$ for some smooth function $\rho$  on $I$.

On the other hand, it follows from 
$$0=d\nu(X)=-\<AX,T\>\,\,\,\mbox{for all}\,\,\, X\in TM^n$$
that $AT=0$, hence ${F}(x,\cdot)\colon\, I\to \Q_\e^n\times \R$ are geodesics in $\Q_\e^n\times \R$, where ${F}=f\circ\psi$. Therefore, the projections $\Pi_1\circ F(x,\cdot)\colon\, I\to \Q_\e^n$ and $\Pi_2\circ F(x,\cdot)\colon\, I\to \R$  are geodesics of $\Q_\e^n$ and $\R$, respectively. 

 That $\Pi_2\circ F(x,\cdot)\colon\, I\to \R$  are geodesics in  $\R$ just means that $\rho(s)=Bs,$
for some constant $B>0$, 
after possibly a translation in the parameter $s$ and changing $s$ by $-s$. Now define $g\colon\,  M^{n-1}\to \Q_\e^n$ by 
$$g(x)=\Pi_1\circ F(x,0).$$
Rescaling the parameter $s$ so that the geodesics $\Pi_1\circ F(x,\cdot)\colon\, I\to \Q_\e^n$ have unit speed, the fact that they are normal to $g$ at $g(x)$ for any $x\in M^{n-1}$ just says  that 
$$\Pi_1\circ F(x,s)=g_s(x),$$
where $g_s$ denotes the parallel hypersurface to $g$ at a distance $s$.\qed

\begin{remark}{\em The proof of Theorem \ref{thm:constantangle} also applies to hypersurfaces of $\R^{n+1}$ whose unit normal vector field makes a constant angle with a fixed direction $\d/\d t$. Namely, writing $\R^{n+1}=\R^n\times \R$, with the second factor being spanned by $\d/\d t$, it shows that any such hypersurface is either an open subset of an affine subspace $\R^n\times \{t_0\}$ for some $t_0\in \R$, an open subset of a product $M^{n-1}\times \R$, where $M^{n-1}$ is a hypersurface of $\R^n$, or it is locally given by (\ref{eq:constantangle}), where $g_s$ is the family of parallel hypersurfaces to some hypersurface $g$ in the first factor $\R^n$, namely, $g_s(x)=g(x)+sN(x)$ for a  unit vector field $N$ to $g$. A proof of this fact for surfaces in $\R^3$ was given in \cite{nistor}.}
\end{remark}

\section[Nonrotational examples in dimension three]{Nonrotational examples in dimension three}

  Here we use the construction of the previous section to produce  a family of nonrotational hypersurfaces of $\Sf^3\times R$ (respectively, $\Hy^3\times \R$) with constant sectional curvature $c$ for any $c\in (0,1)$ (respectively, $c\in (-1,0)$).  
  
  Given a hypersurface $g\colon\, M^{n-1}\to \Q_\e^n$ and the family $g_s \colon\, M^{n-1}\to \Q_\e^n$ of  parallel hypersurfaces to $g$, an easy computation shows that, whenever $\mbox{cot}_\e(s):=C_\e(s)/S_\e(s)$ is not a principal curvature of $g$ at any $x\in M^{n-1}$,  the shape operator $A_s$ of $g_s$ with respect to the unit normal vector field $N_s$ given by (\ref{eq:normalparallel}) is
  \be\label{eq:shape}A_s=(\mbox{cot}_\e s\,I-A)^{-1}(\mbox{cot}_\e s\, A+\e I).
  \ee
  
  Let $g\colon\, M^{2}\to \Q_\e^3$ be a  surface and let 
    $$f\colon\, M^3:=V\times I\subset M^2\times \R\to \Q_\e^3\times \R\subset \Ee^{5}$$  be defined as in the previous section in terms of $g$. The normal space of $f$, as a submanifold of $\Ee^5$, is spanned by the unit normal vector field $\eta$ given by (\ref{eq:eta}) and by the unit normal vector field $\xi(x,s)=g_s(x)$, which is  normal to $ \Q_\e^3\times \R$ at $f(x,s)$. We have 
  $$a \tilde{\nabla}_X\eta=B{g_s}_*A^sX=Bf_*A^sX$$
  and
  $$a \tilde{\nabla}_{\frac{\d}{\d s}}\eta=\e Bg_s=\e B\xi,$$
  hence the principal curvatures of $A^f_\eta$ are 
  $$-\frac{B}{a}k_1^s,\,\,\,-\frac{B}{a}k_2^s\,\,\,\mbox{and}\,\,\, 0,$$
  where $k_1^s$ and $k_2^s$ are the principal curvatures of $g_s$, the principal curvature $0$ corresponding to the principal direction $\d/\d s$. On the other hand,
   $$\tilde{\nabla}_X\xi={g_s}_*X=f_*X$$
  and
  $$ \tilde{\nabla}_{\frac{\d}{\d s}}\xi=N_s=\frac{1}{a^2}f_*{\frac{\d}{\d s}}-\frac{B}{a}\eta.$$
  Thus, the principal curvatures of $A^f_\xi$ are $-1/a^2$ and $-1$, the first being simple with $\d/\d s$ as principal direction, and the second having multiplicity two with $TV$ as eigenbundle. 
  
  Now assume that $M^2=M_0^2$ is flat. Then, the principal curvatures $k_1$ and $k_2$ of $g$ satisfy $k_1k_2=-\e$ everywhere. By (\ref{eq:shape}), the principal curvatures of $g_s$ with respect to  $N_s$ are
  $$k_i^s=  \frac{\mbox{cot}_\e sk_i+\e}{\mbox{cot}_\e s-k_i},\,\,\,1\leq i\leq 2,$$
  hence  $k_1^sk_2^s=-\e$, that is, $g_s$ is also a flat surface. It follows that the sectional curvature of $M^3$ along $TV$ is
  $$(-\frac{B}{a}k_1^s)(-\frac{B}{a}k_2^s)+\e=\frac{\e}{a^2},$$
  which is also the sectional curvature of $M^3$ along any plane spanned by $\d/\d s$ and a vector $X\in TV$. 
  
  \begin{remark}{\em It is easily seen that if the hypersurface $f$ just constructed is  regarded as a submanifold of $\R^5$ for  $\e=1$, then it does not have any umbilical normal direction at any point. Hence it provides a new example of a constant curvature submanifold of $\R^5$ with codimension two that is free of weak-umbilic points in the sense of \cite{moore}.}
  \end{remark}
  
 \begin{example} {\em As an explicit example, consider the Clifford torus $$g\colon\,M^2_0:=\Sf^1(\cos\theta_0)\times  \Sf^1(\sin\theta_0)\to \Sf^3$$ parametrized by
 $$g(t_1,t_2)=(\cos\theta_0\cos t_1, \cos\theta_0\sin t_1,\sin\theta_0\cos t_2, \sin\theta_0\sin t_2),$$
 which has
 $$N(t_1,t_2)=(-\sin\theta_0\cos t_1, -\sin\theta_0\sin t_1,\cos\theta_0\cos t_2, \cos\theta_0\sin t_2)$$
as a unit normal vector field in $\Sf^3$. Then,
$$f\colon\,M^2_0\times \R\to \Sf^3$$
given by (\ref{eq:constantangle}) can be reparametrized by
$$ f(t_1,t_2,s)=(\cos s\cos t_1, \cos s\sin t_1,\sin s\cos t_2, \sin s\sin t_2,Bs),$$
after replacing $s+\theta_0$ by $s$ and a translation in the $\d/\d t$-direction. This hypersurface appears in \cite{dt} as an example of a weak-umbilic free  doubly-rotation surface with constant sectional curvature having the helix $s\mapsto (\cos s, \sin s, Bs)$ as profile, in the sense of \cite{dillen2}. 

A similar example can be constructed in $\Hy^3\times \R$, starting with the flat surface $$g\colon\,M^2_0:=\Hy^1(\cosh\theta_0)\times  \Sf^1(\sinh\theta_0)\to \Hy^3$$ parametrized by
 $$g(t_1,t_2)=(\cosh\theta_0\cos t_1, \cosh\theta_0\sin t_1,\sinh\theta_0\cos t_2, \sinh\theta_0\sin t_2).$$
 In this case, the corresponding constant curvature hypersurface of   $\Hy^3\times \R$ is
 $$ f(t_1,t_2,s)=(\cosh s\cos t_1, \cosh s\sin t_1,\sinh s\cos t_2, \sinh s\sin t_2,Bs),$$
These examples can be characterized as the only  constant curvature hypersurfaces of   $\Q_\e^3\times \R$ with $0$ as principal curvature in the $T$-direction and whose two remaining principal curvatures are constant along $\{T\}^\perp$. }
 \end{example}

\section[The main result]{The main result}

In this section we prove our main result, namely, we provide a complete classification of all hypersurfaces with constant sectional curvature of 
$\Q_\e^n\times \R$, $n\geq 3$. We state separately the cases $\e=1$ and $\e=-1$. For $\e=1$ we have:
\begin{theorem}\po\label{thm:sphere} Let $f\colon\, M_c^n\to \Sf^n\times \R$, $n\geq 3$, be an isometric immersion of a Riemannian manifold of constant sectional curvature $c$. Then $c\geq 0$. Moreover,
\begin{itemize}
\item[$(i)$] if $c=0$ then $n=3$ and $f(M^3_0)$ is an open subset of a Riemannian product $M^2_0\times \R$, where $M^2_0$ is a flat surface of $\Sf^3$.
\item[$(ii)$] if $c\in (0,1)$ then $n=3$ and $f$ is locally given by the construction described in  Section $6$. 
\item[$(iii)$] if $c=1$ then $f(M^n_1)$ is an open subset of a slice $\Sf^n\times \{t\}$. 
\item[$(iv)$] if $c>1$ then $f(M^n_c)$ is an open subset of a rotation hypersurface given by Theorem \ref{thm:rotationcsc1}-$(ii)$.
\end{itemize}
\end{theorem}

The classification of  constant curvature hypersurfaces of $\Hy^n\times \R$ with dimension $n\geq 3$ reads as follows:

\begin{theorem}\po\label{thm:hyperbolic} Let $f\colon\, M_c^n\to \Hy^n\times \R$, $n\geq 3$, be an isometric immersion of a Riemannian manifold of constant sectional curvature $c$. Then $c\geq -1$. Moreover,
\begin{itemize}
\item[$(i)$] if $c=-1$ then $f(M^n_{-1})$ is an open subset of a slice $\Hy^n\times \{t\}$.
\item[$(ii)$] if $c\in (-1,0)$ then either $n=3$ and $f$ is locally given by the construction described in  Section $6$, or $f(M^n_0)$ is an open subset of one of the  rotation hypersurfaces given by Theorem \ref{thm:rotationcsc2}-$(ii)$.
\item[$(iii)$] if $c=0$ then one of the following possibilities holds:
\begin{itemize}
\item[$(a)$] $n=3$ and $f(M^3_0)$ is an open subset of a Riemannian product $M^2_0\times \R$, where $M^2_0$ is a flat surface of $\Hy^3$.
\item[$(b)$]  $f(M^n_0)$ is an open subset of a Riemannian product $M^{n-1}_0\times \R$, where $M^{n-1}_0$ is a horosphere of $\Hy^n$.
\item[$(c)$] $f(M^n_0)$ is an open subset of the spherical rotation hypersurface given by Theorem \ref{thm:rotationcsc2}-$(iii)$-$(a)$.
\end{itemize} 
\item[$(iv)$] if $c>0$ then $f(M^n_c)$ is an open subset of the spherical rotation hypersurface given by Theorem \ref{thm:rotationcsc2}-$(iv)$.
\end{itemize}
\end{theorem}

\noindent{\em Proof of Theorems \ref{thm:sphere} and \ref{thm:hyperbolic}:} Assume that the vector field $T$ does not vanish at $x\in M^n$. Then $T$ is a principal direction of $f$ by 
Lemma \ref{le:Tprincipal1} and Lemma \ref{le:Tprincipal2}. 
Let $\{T,X_1,\ldots,X_{n-1}\}$ be an orthogonal basis of
eigenvectors of $A_N$ at $x$, with $$A_NT=\lambda T\,\,\,\mbox{ and}\,\,\,A_NX_i=\lambda_iX_i,\,\,\,1\leq i\leq n-1.$$   From the Gauss equation (\ref{eq:Gauss}) of $f$ for $X=X_i$ and $Y=Z=X_j$, $i\neq j$, we get
\be\label{eq:gauss1} c-\e =\lambda_i\lambda_j,\,\,\,\,\,i\neq j.\ee
On the other hand, for $X=T$ and $Y=Z=X_i$ the 
Gauss equation yields
\be\label{eq:gauss2}
c-\e=\lambda\lambda_i-\e||T||^2.
\ee
Assume first that $c=\e$. By (\ref{eq:gauss1}), we can assume that $\lambda_i=0$ for all $2\leq
i\leq n-1$. Then, applying (\ref{eq:gauss2}) for $i\geq 2$ yields a contradiction with $T\neq 0$.  We conclude that for $c=\e$ the vector field $T$ vanishes identically, and this gives part $(iii)$ of Theorem \ref{thm:sphere} and part $(i)$ of Theorem \ref{thm:hyperbolic}.

Now suppose that $c\neq \e$. Then  $T$ can not vanish on any open subset. Thus,  we can assume without loss of generality that it is nowhere vanishing. 
If $n\geq 4$, we obtain from (\ref{eq:gauss1}) that all $\lambda_i's$ coincide for $2\leq i\leq n-1$. Denote all of  them by $\mu$. Then,  the Gauss equations now read
\be\label{eq:G1}
c-\e=\mu^2
\ee
and
\be \label{eq:G2}
c-\e=\lambda\mu-\e\|T\|^2,
\ee
which can also be written as
\be \label{eq:G2b}
c=\lambda\mu+\e\nu^2.
\ee
In particular, it follows from (\ref{eq:G1}) that $c>\e$.

Now, since $T\neq 0$, it follows from (\ref{eq:G1}) and  (\ref{eq:G2}) that $\lambda\neq \mu$. Moreover, since $T$ is a principal direction, we obtain  from (\ref{eq:DerivadaNu})  that $\nu$ is constant along the leaves of $\{T\}^\perp$, and hence the same holds for $\lambda$ by (\ref{eq:G2b}) (since $\mu$ has multiplicity greater than one, one can show using the Codazzi equation that it is constant along its eigenbundle; cf. the proof of Theorem 1 in \cite{veken}). Then, one can use the following result to conclude that $f$ is a rotation hypersurface. It slightly generalizes   Theorem 1 in  \cite{veken}, but actually follows from its proof. 

\begin{proposition}\label{prop:rotation}
Let $f:M^n\to\Q_\e^n\times\R$ be a hypersurface with $n\geq 3$ and $T\neq 0$. Assume that $f$ has exactly two principal curvatures $\lambda$ and $\mu$ everywhere, the first one being simple with $T$ as a principal direction. If $\lambda$ is constant along the leaves of the eigenbundle $\{T\}^\perp$ of $\mu$, then $f(M^n)$ is an open subset of a rotation hypersurface.
\end{proposition}

Thus, the proofs of Theorems \ref{thm:sphere} and \ref{thm:hyperbolic} for $n\geq 4$ are completed by Theorems \ref{thm:rotationcsc1} and \ref{thm:rotationcsc2}. This also applies to the case  $n=3$ when we have $\lambda_2=\lambda_3$ everywhere. By  (\ref{eq:gauss1}) and  (\ref{eq:gauss2}), this is not the case only if   $\lambda=0$. In this situation, equation (\ref{eq:G2b}) reduces to 
\be\label{eq:G2c}\e\nu^2=c.\ee 
If $c=0$, then $\nu$ vanishes identically, and thus $f(M_0^3)$ must be an open subset of a Riemannian product $M^2_0\times \R$, where $M^2_0$ is a flat surface in either $\Sf^3$ or $\Hy^3$, according as $\e=1$ or $\e=-1$, respectively. If $c\neq 0$, it follows from (\ref{eq:G2c}) that $f$ is a constant angle hypersurface. Therefore, by Theorem \ref{thm:constantangle} it is locally given by (\ref{eq:constantangle}) for some surface $g\colon\, M^{2}\to \Q_\e^3$. Moreover, if we write $\nu=1/a$, it was shown in Section $6$ that the principal curvatures of $f$ are 
$$-\frac{B}{a}k_1^s\,-\frac{B}{a}k_2^s\,\,\,\mbox{and}\,\,\, 0,$$
  where $k_1^s$ and $k_2^s$ are the principal curvatures of $g_s$. By the Gauss equation (\ref{eq:gauss1}), we have
  $$c-\e=(-\frac{B}{a}k_1^s)(-\frac{B}{a}k_2^s).$$
  Replacing $c=\e/a^2$ and using that $B^2+1=a^2$, it follows that 
  $k^s_1k^s_2=-\e$, hence  $g$ is a flat surface.\qed

\end{document}